\documentclass[twoside,reqno]{amsart} 
\usepackage{amsthm}           
\usepackage{xspace}
\usepackage[matrix,arrow]{xy}          
\xymatrixcolsep{1.7pc} 
\xymatrixrowsep{1.7pc}
\newdir{ >}{{}*!/-5pt/\dir{>}}  
\newdir{> }{{}*!/-9pt/\dir{>}}  

\newcommand{\ie}                {\emph{i.e.,}\xspace}

\newcommand{\tensor}            {\otimes}
\newcommand{\Smash}             {\wedge}

\newcommand{\Wedge}             {\vee}

\newcommand{\ra}        {\rightarrow}

\newcommand{\hra}       {\hookrightarrow}

\DeclareMathOperator{\Hom}{Hom}
\DeclareMathOperator{\Ext}{Ext}
\DeclareMathOperator{\PExt}{PExt}

\newcommand{\field}[1]  {\mathbb{#1}}

\newcommand{\Z}         {\field Z}
\newcommand{\Q}         {\field Q}
\newcommand{\F}         {\field F}

\newcommand{\FF}{\textup{F}}

\newcommand{\donothing}[1]{}

\numberwithin{equation}{section}

\usepackage{enumerate}        
\newenvironment{roenumerate}{\begin{enumerate}[\upshape (i)]}{\end{enumerate}}
\newcommand{\ignore}[1]{}

\newtheorem{thm}[equation]{Theorem} 
\newtheorem{defn}[equation]{Definition}
\newtheorem{prop}[equation]{Proposition}
\newtheorem{cor}[equation]{Corollary}
\newtheorem{lemma}[equation]{Lemma}

\theoremstyle{definition}

\newcommand{\ulp}{\textup{(}}
\newcommand{\urp}{\textup{)}}
\newcommand{\uc}{\textup{:}}

\newcommand{\dfn}{\textbf} 
\newcommand{\mdfn}[1]{\dfn{\mathversion{bold}#1}} 

\newcommand{\mathcolon}{\colon\,} 
\newcommand{\mc}{\mathcolon}

\newcommand{\EM}                {Eilenberg--Mac\,Lane\xspace}

\newcommand{\Zn}{\Z \mspace{-1.0mu}/\mspace{-1.0mu} n}

\newcommand{\Zpl}{\Z \mspace{-1.0mu}/\mspace{-1.0mu} p^{l}}

\newcommand{\Fp}{\F_p}

\newcommand{\Hp}{H\Fp}

\newcommand{\Ph}                {\textup{Ph}}

\newcommand{\oPh}{\textup{$1$-Ph}}
\newcommand{\nPh}{\textup{$n$-Ph}}
\newcommand{\cat}[1]{\mathcal{#1}}
\newcommand{\cD}{\cat{D}}

\begin{document}

\title{Phantom maps and chromatic phantom maps}

\author{J. Daniel Christensen}
\address{Department of Mathematics\\ Johns Hopkins University\\ 
Baltimore, MD 21218} 
\email{jdc@math.jhu.edu}
\author{Mark Hovey}
\address{Department of Mathematics\\ Wesleyan University\\ 
Middletown, CT 06459} 
\email{hovey@member.ams.org}

\subjclass{Primary 55P42}

\keywords{Phantom map, chromatic phantom map, $n$-phantom map,
cohomotopy, stable homotopy, spectrum, $n$-finite type}

\begin{abstract}
In the first part, we determine conditions on spectra $X$ and $Y$ 
under which either every map from $X$ to $Y$ is phantom, or no
nonzero maps are.
We also address the question of whether such all or nothing behaviour
is preserved when $X$ is replaced with $V \Smash X$ for $V$ finite.
In the second part, we introduce chromatic phantom maps.  A map is
$n$-phantom if it is 
null when restricted to finite spectra of type at least $n$.
We define divisibility and finite type conditions which are
suitable for studying $n$-phantom maps.  We show that the
duality functor $W_{n-1}$ defined by Mahowald and Rezk is the analog
of Brown-Comenetz duality for chromatic phantom maps, and give conditions under which the
natural map $Y \ra W_{n-1}^{2} Y$ is an isomorphism.
\end{abstract}

\date{\today}

\maketitle

\tableofcontents

\renewcommand{\baselinestretch}{1.2}\normalsize

\section{Introduction}\label{se:intro}

This paper takes place in the stable homotopy category.
A map $X \ra Y$ of spectra is said to be \dfn{phantom} if for every
finite spectrum $V$ and every map $V \ra X$ the composite $V \ra X \ra Y$
is zero.
We write $\Ph(X,Y)$ for the group of all phantom maps from $X$ to $Y$.

The first half of the paper, Section~\ref{se:all}, deals with some
results of a general nature.
The theme is to determine conditions on spectra $X$ and $Y$ 
under which either every map from $X$ to $Y$ is phantom, or no
nonzero maps are.
We also address the question of whether such all or nothing behaviour
is preserved when $X$ is replaced with $V \Smash X$ for $V$ finite.
Our main tool is the observation that for $Y$ of finite type,
a map $X \ra Y$ is phantom if and only if it is a divisible
element of the group $[X,Y]$.
We also make use of the connection between phantom maps and
Brown-Comenetz duality.

In the second half of the paper, Section~\ref{se:n-ph}, 
we work in the $p$-local stable homotopy category and introduce 
generalizations of phantom maps called chromatic phantom maps.
A map $X \ra Y$ is \mdfn{$n$-phantom} if for every 
finite spectrum $V$ of type at least $n$ and every map
$V \ra X$ the composite $V \ra X \ra Y$ is zero.
The results of this section, which require more background of the reader, 
generalize some of the results of the first half of the paper.
We define analogs of the divisibility and finite type conditions 
which are suitable for studying chromatic phantom maps.
The significance of chromatic phantom maps is demonstrated by the fact 
that the kernel in homotopy of the important map
$L_{n}S^{0}\ra L_{K(n)}S^{0}$ is precisely the $n$-phantom
homotopy classes.  
We show that the duality functor $W_{n-1}$ defined by Mahowald and
Rezk~\cite{mare} is the analog of Brown-Comenetz duality.
Our Theorem~\ref{th:n-ph} gives conditions under which the
natural map $Y \ra W_{n-1}^{2} Y$ is an isomorphism, generalizing
a result of Mahowald and Rezk.

Special thanks are due to Neil Strickland for
asking the question which led to Section~\ref{ss:nothing}
and to Sharon Hollander for asking the question which 
led to many of the results in the first half of the paper.
We also appreciate helpful conversations with Amnon Neeman.
Correspondence with Haynes Miller led to the results in
Section~\ref{se:all}, and we thank him for sharing his
insights, ideas and arguments.

\section{Phantom maps}\label{se:all}

\subsection{Divisibility}

In this section we describe conditions on spectra $X$ and $Y$
which ensure that all maps from $X$ to $Y$ are phantom.  
This will turn out to have close connections to the divisibility
of the group $[X,Y]$.

We write $IY$ for the Brown-Comenetz dual of $Y$~\cite{brco:pdg}, 
which is characterized by the natural
isomorphism $[X,IY] = \Hom(\pi_{0}(X \Smash Y),\Q/\Z)$.  
We write $I$ for $IS^{0}$, and note that $IY=\FF(Y,I)$, where 
$\FF(-,-)$ denotes the function spectrum.
$I(-)$ is a contravariant functor on the stable homotopy category, 
and there is a natural map $Y \ra I^{2}Y$.
As a special case of the defining property, we have 
$\pi_{k}IY = \Hom(\pi_{-k}Y,\Q/\Z)$.

We will make use of the following proposition which gives several
characterizations of phantom maps.
This is Proposition~4.12 from~\cite{chst:pmht}.

\begin{prop}\label{pr:char}
The following conditions on a map $X \ra Y$ are equivalent\uc 
\begin{roenumerate}
\item For each finite spectrum $V$ and each map $V \ra X$, the
composite $V \ra X \ra Y$ is zero \ulp\ie $X \ra Y$ is phantom\urp.
\item For each homology theory $h$ taking values in abelian groups, 
the map $h(X) \ra h(Y)$ is zero.
\item\label{it:char} The composite $X \ra Y \ra I^{2}Y$ is zero.   \qed
\end{roenumerate}
\end{prop}

Part~(iii) of Proposition~\ref{pr:char} implies that the fiber $F\ra Y$
of the map $Y\ra I^{2}Y$ is phantom.  Hence, there are no nonzero
phantom maps to $Y$ if and only if the map $Y\ra I^{2}Y$ is a split
monomorphism.  This observation gives us some information about the
homotopy groups of such a spectrum $Y$.  Indeed, if there are no nonzero
phantom maps to 
$Y$, then each homotopy group $\pi _{n}Y$ must split off its double dual
$\Hom (\Hom (\pi _{n}Y,\Q /\Z),\Q /\Z)$.  Furthermore, each summand of
$\pi _{n}Y$ must split off its double dual.  In particular, if $\Z $ is
a summand of $\pi _{n}Y$, there must be a nonzero phantom map to $Y$.
Indeed, if not, $\Z $ would be a summand of its double dual $\hat{\Z}$,
the pro-finite completion of $\Z$.  Denote the other summand of
$\hat{\Z}$ by $A$.  Then we would have $A/p=0$ for all primes $p$,
and so $A$ would have to be a divisible group.  But there are no
divisible summands of $\hat{\Z}$, so $A=0$, and we find that $\Z\cong
\hat{\Z}$.  This is a contradiction.  

Recall that a spectrum $Y$ is of \dfn{finite type} if each homotopy
group $\pi _{n}Y$ is a finitely generated abelian group.  
A finite type spectrum $Y$ is rationally trivial if and only if
its homotopy groups are finite, \ie if and only if the map $Y \ra I^{2}Y$
is an isomorphism.
Thus for $Y$ of finite type, there are no nonzero phantom maps to $Y$ 
if and only if $Y$ is rationally trivial.  

Next we present a new characterization of phantom maps,  
based on this observation.  
A map $f \mathcolon X \ra Y$ is \dfn{divisible} if it is a divisible element
of the group $[X,Y]$, \ie if for each nonzero integer $n$,
we have $f = n g$ for some $g$.

\begin{thm}\label{th:div}
Let $Y$ be a spectrum of finite type.  Then a map $X \ra Y$ is phantom
if and only if it is divisible.
In particular, all maps from $X$ to $Y$ are phantom if and only if
$[X,Y]$ is divisible.  
Moreover, the group $\Ph(X,Y)$ of phantom maps is divisible.
\end{thm}

\begin{proof}
Let $f \mathcolon X \ra Y$ be a divisible map and let $V \ra X$ be a map from
a finite.  Then the composite $V \ra X \ra Y$ is a divisible
element of the group $[V,Y]$.  This group is finitely generated
(since $Y$ has finite type), and so the composite is zero.
Thus $f$ is phantom.

Now suppose that $f \mathcolon X \ra Y$ is phantom.
By part~(iii) of Proposition~\ref{pr:char}, the given map $f$
factors through the fibre $F$ of the natural map $Y \ra I^{2}Y$.
Since $Y$ has finite type, the fibre $F$ is 
rational~\cite[Lemma~A3.13]{ma:ssa}.
Thus $f$ must be divisible.

In fact, for any nonzero integer $n$, we can write $f$ as $n g$,
with $g$ \emph{phantom}, since the map $F \ra Y$ is phantom (by
Proposition~\ref{pr:char}).
This proves that $\Ph(X,Y)$ is divisible.

We have proved the first and third statements.  From the first it is
clear that all maps from $X$ to $Y$ are phantom if and only if $[X,Y]$
is divisible.
\end{proof}

The proof shows that any phantom map $X \ra Y$,
with $Y$ of finite type, factors through a rational spectrum,
and therefore through the rationalization of $X$.
This is a stable analogue of~\cite[Theorem~5.1~(ii)]{mc:pm}.

Phantom maps are not always divisible.
For example, there is a spectrum $Y$ such that $\Ph(\Hp,Y)$
is nonzero~\cite[Proposition~4.18]{chst:pmht}.
Of course, by the above, no such $Y$ has finite type.
Also, a divisible map is not always phantom.
For example, the identity map of $H\Q$ is divisible but not phantom.

A natural question is the following.  Suppose that all graded maps
from $X$ to $Y$ are phantom, \ie that $\Ph(X,Y)_{*} = [X,Y]_{*}$.
Does it follow that all maps from $V \Smash X$ to $Y$ are phantom
for $V$ finite?
The answer to this question is no, as we will see after the following
theorem. 

\begin{thm}\label{th:rat}
Let $Y$ be a spectrum of finite type.  Then the following are
equivalent.
\begin{roenumerate}
\item $\Ph(V \Smash X,Y)_{*} = [V \Smash X,Y]_{*}$ for all finite $V$.
\item $\Ph(V \Smash X,Y)_{*} = [V \Smash X,Y]_{*}$ for all $V$.
\item $X \Smash IY = 0$.
\item $[X,Y]_{*}$ is rational.
\end{roenumerate}
\end{thm}

\begin{proof}
(iii) $\implies$ (ii):  Assume that $X \Smash IY = 0$.  Then
$V \Smash X \Smash IY = 0$ for any $V$.
So $[V \Smash X,I^{2}Y]_{*} = \Hom(\pi_{-*}(V \Smash X \Smash IY),\Q/\Z) = 0$.
So every map from $V \Smash X$ to $Y$ is phantom, by part~(\ref{it:char}) of
Proposition~\ref{pr:char}.

(ii) $\implies$ (i):  This is clear.

(i) $\implies$ (iv):  Assume that all graded maps $V \Smash X \ra Y$
are phantom for finite $V$.
Let $p$ be a prime and write $M(p)$ for the mod $p$ Moore spectrum.
By assumption, any graded map $M(p) \Smash X \ra Y$ is phantom.
By Theorem~\ref{th:div}, any such map is divisible.  But since $M(p)$
is torsion, we conclude that any map $M(p) \Smash X \ra Y$ is zero.
Therefore, multiplication by $p$ on $[X,Y]_{*}$ is an isomorphism.
This is true for all primes, so $[X,Y]_{*}$ is rational.

(iv) $\implies$ (iii):  
Assume that $[X,Y]_{*}$ is rational.
By the following lemma,
$\FF(X,I^{2}Y)$ is the
pro-finite completion of $\FF(X,Y)$. 
By hypothesis, $\FF(X,Y)$ is rational, and so its pro-finite completion is
trivial. 
Thus $\Hom(\pi_{*}(X \Smash IY),\Q/\Z) = 0$.
Since $\Q/\Z$ is a cogenerator, it follows that $\pi_{*}(X \Smash IY) = 0$.
Therefore $X \Smash IY = 0$.
\end{proof}

We write $Z \ra \hat{Z}$ for the pro-finite completion of a spectrum $Z$.
This is Bousfield localization with respect to the wedge over all
primes of the mod $p$ Moore spectra.
Put another way, $Z \ra \hat{Z}$ is the initial map to a spectrum
which has no nonzero maps from a rational spectrum.
We write $\hat{A}$ for the pro-finite completion of an abelian group $A$.

The above proof used the following lemma.

\begin{lemma}\label{le:prof}
For $Y$ finite type and any $X$, the natural map 
$\FF(X,Y) \ra \FF(X,I^{2}Y)$ is pro-finite completion.
In particular, $Y \ra I^{2}Y$ is pro-finite completion.
\end{lemma}

\begin{proof}
The second statement is~\cite[Theorem~9.11]{ma:ssa}.  The general
case follows.  Indeed, for a general spectrum $E$, if $Y$ is $E$-local,
so is $\FF(X,Y)$ for any $X$.  In particular, $\FF(X,I^{2}Y)$ is local with
respect to the wedge of all the $M(p)$.  Since the fiber $F$ of
$Y\xrightarrow{}I^{2}Y$ is rational, so is $\FF(X,F)$.  Thus $\FF(X,F)$ is
acyclic with respect to the wedge of all the $M(p)$.  
\end{proof}

We now find specific spectra $X$ and $Y$ where $Y$ is of finite type and
$[X,Y]_{*}$ is divisible, but not rational.  In view of
Theorems~\ref{th:div} and~\ref{th:rat}, this implies that all maps from 
$X$ to $Y$ are
phantom, but that there is a map $M(p)\Smash X\xrightarrow{}Y$ which is
not phantom.

\begin{prop}\label{pr:counter}
The integral cohomology of $\hat{S^{0}}$ is $\Q /\Z $, concentrated in
degree $1$.  In particular, all maps from $\hat{S^{0}}$ to $H\Z $ are
phantom, but there are non-phantom maps from $\hat{S^{0}}\Smash
M(p)=M(p)$ to $H\Z$.
\end{prop}

\begin{proof}
We use the universal coefficient theorem to calculate
$H^{*}(\hat{S^{0}})$ from the fact that $H_{*}(\hat{S^{0}})=\hat{\Z}$,
the pro-finite completion of $\Z $, concentrated in degree $0$.  We have
isomorphisms $H^{0}(\hat{S^{0}})\cong \Hom (\hat{\Z},\Z)$ and
$H^{1}(\hat{S^{0}})\cong \Ext (\hat{\Z},\Z )$, while the other groups
vanish.  Suppose we have a
nonzero homomorphism $f\mathcolon \hat{\Z}\ra \Z$.  Then the
image of $f$ must be the subgroup $m\Z$ for some nonzero $m$.  If we
define $g=(\frac{1}{m})f$, then $g$ is a well-defined homomorphism whose
image is all of $\Z$.  By choosing a section of $g$, we find that
$\hat{\Z}\cong \Z \oplus A$ for some abelian group $A$.  We have already
seen that this is impossible in the paragraph following
Proposition~\ref{pr:char}.  Hence $\Hom (\hat{\Z},\Z)=0$.

In order to calculate $H^{1}(\hat{S^{0}})$, we use some topology.
First, $(H\hat{\Z})^{*}(\hat{S^{0}})\cong (H\hat{\Z})^{*}(S^{0})\cong
\hat{\Z}$ concentrated in degree $0$.  
In particular, $(H\hat{\Z})^{1}(\hat{S^{0}})$ is zero.
It follows from Proposition~\ref{pr:char}~(iii) that the 
group $H^{1}(\hat{S^{0}})$ is all
phantom, and so is divisible.  Moreover, any non-torsion
element of $H^{1}(\hat{S^{0}})$ would survive to give a nonzero element
of $(H\Q )^{1}(\hat{S^{0}})=0$.  Thus $H^{1}(\hat{S^{0}})$ is a direct
sum of copies of $\Q /\Z _{(p)}$ for various primes $p$.  On the other
hand, $(H\Fp )^{*}(\hat{S^{0}})=\Fp $ concentrated in degree $0$.  Hence
there must be exactly one summand of the form $\Q /\Z _{(p)}$ for each
prime $p$, and so $H^{1}(\hat{S^{0}})=\Q /\Z $, as claimed.  
\end{proof}

There are certain situations where the phenomenon of
Proposition~\ref{pr:counter} can not occur. 
We say that a spectrum is \dfn{bounded above} if its
homotopy groups are bounded above.

\begin{thm}\label{th:bdd-above}
Let $Y$ be a spectrum.
Then $[\Hp,Y]^{*} = 0$ for all primes $p$ if and only if 
$[X,Y]$ is rational for each bounded above $X$.
Moreover, if $Y$ has finite type, then a third equivalent condition is 
that all maps from every bounded above $X$ to $Y$ are phantom.
\end{thm}

\begin{proof}
It is clear that $[\Hp,Y]^{*} = 0$ for all primes $p$ if and only if
$[H\Z,Y]$ is rational.  
If $[H\Z,Y]$ is rational, then so is $[X,Y]$ for any $X$ in the
localizing subcategory generated by $H\Z$, since rational abelian
groups are closed under kernels, cokernels, extensions and products.
If $X$ is bounded above then $X$ is the homotopy colimit of its
Postnikov sections, each of which has only finitely many nonzero
homotopy groups.  A spectrum with only finitely many nonzero homotopy
groups is in the thick subcategory generated by $H\Z$, so a bounded
above spectrum is in the localizing subcategory generated by $H\Z$.
This proves the first statement.

Now assume that $Y$ has finite type.
If $[X,Y]$ is rational for bounded above $X$, then by Theorem~\ref{th:div} 
all maps from every bounded above $X$ to $Y$ are phantom.
And if all maps from every bounded above $X$ to $Y$ are phantom,
then using Theorem~\ref{th:div} again we see that $[\Hp,Y]$ is divisible 
and hence zero.  This proves the second statement.
\end{proof}

Note that our proof of
Theorem~\ref{th:bdd-above} applies more generally to spectra $X$ in the
localizing subcategory generated by $H\Z$.  

When $Y$ is harmonic, $[\Hp,Y]^{*} = 0$ for all 
$p$~\cite[Theorem~4.8]{ra:lrc}, and this gives rise to many examples.
Recall that Hopkins and Ravenel proved that suspension
spectra are harmonic~\cite{hora:ssh}.  
Therefore, a finite type suspension spectrum
or a desuspension of a finite type suspension spectrum
satisfy the hypotheses on $Y$ in the theorem.
In particular, a finite spectrum satisfies the hypotheses on $Y$
and so a special case of Theorem~\ref{th:bdd-above} is that all maps from
an \EM spectrum to a finite spectrum are phantom.

We now deduce a variant.  Note that an \EM spectrum has bounded below
cohomotopy groups~\cite[Theorem~4.2]{li:dems}.

\begin{cor}\label{co:cohom}
Let $X$ have bounded below cohomotopy groups and let $Y$ be a 
finite spectrum.  Then all maps from $X$ to $Y$ are phantom,
and $[X,Y] = \Ph(X,Y)$ is rational.
\end{cor}

In particular, if $X$ has bounded below cohomotopy groups,
then these groups consist of phantom maps and are rational.

\begin{proof}
Let $V$ be a finite spectrum and let $V \ra X$ be a map.  
We must show that the restriction map $[X,Y] \ra [V,Y]$ is zero,
so it suffices to prove that the map $\FF(X,Y) \ra \FF(V,Y)$
is zero in homotopy groups.
But, by induction over the cells of $Y$, one can prove that $\FF(X,Y)$
has bounded above homotopy groups.
And since $\FF(V,Y)$ is finite, the result follows from 
Theorem~\ref{th:bdd-above}.
\end{proof}

\subsection{Rational homotopy}

In this section we give stable analogs of results of McGibbon and Roitberg.
First we have the analog of part~(i) of~\cite[Theorem~2]{mcro:pmre}.

\begin{prop}\label{pr:monic}
If $X \ra X'$ induces a monomorphism on rational homotopy, and if
$Y$ is of finite type, then the natural map $\Ph(X',Y) \ra \Ph(X,Y)$
is surjective.
\end{prop}

\begin{proof}
Let $f \mathcolon X \ra Y$ be a phantom map.
By the proof of Theorem~\ref{th:div}, $f$ factors through the 
rationalization $X_{\Q}$ of $X$.
Since $f$ induces a monomorphism on rational homotopy groups, 
the map $f_{\Q} \mathcolon X_{\Q} \ra X'_{\Q}$ is split monic.
Thus $f$ factors through the composite $X \ra X_{\Q} \ra X'_{\Q}$.
And so $f$ extends over $X'$ by composing with the universal map 
$X' \ra X'_{\Q}$.
The resulting map $X' \ra X'_{\Q} \ra Y$ is divisible and hence phantom.
\end{proof}

As a consequence we get a stable version of part 
of~\cite[Theorem~1]{mcro:pmre}.

\begin{cor}\label{co:rat}
If there exists a rational monomorphism from $X$ to a wedge of 
finite spectra, then $\Ph(X,Y) = 0$ for all finite type $Y$. \qed
\end{cor}

Note that if $X \ra X'$ induces an \emph{epimorphism} on rational
homotopy, we can not conclude that $\Ph (X',Y)\xrightarrow{}\Ph (X,Y)$
is a monomorphism.  Indeed, consider the rational isomorphism
$S^{0}\xrightarrow{}H\Q $.  Then $\Ph(S^{0},Y)=0$ for any $Y$, but
$\Ph(H\Q ,Y)$ is often nonzero, for example when $Y=S^{1}$.  

Next we have the stable analog of part~(ii) of~\cite[Theorem~2]{mcro:pmre}.

\begin{prop}\label{pr:epi}
If $Y$ and $Y'$ are of finite type and $f\mc Y\ra Y'$ induces an epimorphism
on rational homotopy, then the natural map $\Ph (X,Y)\ra \Ph (X,Y')$ is
surjective for all $X$.  
\end{prop}

\begin{proof}
Let $F$ denote the fiber of $Y\ra I^{2}Y$ and let $F'$ denote the fiber
of $Y'\ra I^{2}Y'$.  Then there is a unique induced map $F\ra F'$, since
$F$ is acyclic and $I^{2}Y'$ is local with respect to the wedge of all
the $M(p)$ as $p$ runs through the primes.  We first claim that the map
$I^{2}f \mathcolon I^{2}Y \ra I^{2}Y'$ is also an epimorphism on rational homotopy.
Indeed, by hypothesis, the cokernel of $\pi _{n}f$ is a (necessarily
finite) torsion group.  Since applying $I^{2}$ tensors the homotopy
of a finite type spectrum with $\hat{\Z}$, the cokernel of 
$\pi_{n}I^{2}f$ is also a bounded torsion group, and so $I^{2}f$ induces an
epimorphism on rational homotopy.  The five lemma then implies that $F
\ra F'$ induces an epimorphism on rational homotopy, since we have a
short exact sequence 
\[
0 \ra \pi _{n}Y \ra \pi _{n}I^{2}Y \ra \pi _{n-1} F \ra 0
\]
and a similar one for $Y'$.  Since $F$ and $F'$ are rational spectra, the
map $F \ra F'$ is a split epimorphism.  Given a phantom map $X \ra Y'$,
it must factor through a map $X \ra F'$.  This map then factors through
$F$, and so the original phantom map factors through $Y$, as required.  
\end{proof}

Once again, if $Y \ra Y'$ induces a \emph{monomorphism} on rational
homotopy, we can not conclude that $\Ph (X,Y) \ra \Ph(X,Y')$ is a
monomorphism.  Take $Y=Y'=H\Z $, with the map $Y \ra Y'$ being
multiplication by $n$, which is a rational isomorphism.  Take
$X=\hat{S^{0}}$.  Then we have seen in Proposition~\ref{pr:counter} that
$\Ph (X,Y)=\Q /\Z $, on which multiplication by $n$ is surjective but not
injective.

\subsection{Preserving nothing}\label{ss:nothing}

The final part of this section is concerned with the following question, 
which is due to Neil Strickland.
Suppose we know that no maps of any degree from $X$ to $Y$ are
phantom.  Is it necessarily the case that no maps from $V \Smash X$
to $Y$ are phantom for any finite $V$?

We know of no case where the answer is no.

The answer is yes when $Y$ is an \EM spectrum
and also when $Y$ is a finite type spectrum.  We begin with the \EM spectrum
case. 

\begin{thm}\label{th:PExt}
Let $X$ be a spectrum and $G$ an abelian group such that
there are no phantom maps of any degree from $X$ to $HG$.
Then there are no phantom maps from $V \Smash X$ to $HG$ for any finite $V$.
\end{thm}

In order to prove this, we need to understand phantom maps from
a spectrum $X$ to an \EM spectrum.
Recall that the universal coefficient theorem gives us a (split) 
monomorphism $\Ext(H_{-1}X,G) \hra [X,HG]$ whose image consists 
of those maps sent to zero by the integral homology functor $H_{0}(-)$.
Thus the group of phantom maps from $X$ to $HG$ can be thought of
as a subgroup of $\Ext(H_{-1}X,G)$.
Theorem~6.4 of~\cite{chst:pmht} identifies this subgroup as
$\PExt(H_{-1}X,G)$, where we write $\PExt(A,B)$ for the subgroup of 
$\Ext(A,B)$ consisting of the divisible elements.
(These correspond to the ``pure'' extensions;  
see~\cite[Section~6]{chst:pmht}.)
Since the monomorphism above is split, divisibility in the $\Ext$
group is equivalent to divisibility in $[X,HG]$, so we can make the 
slightly sharper statement that a map $X \ra HG$ is phantom if and 
only if it is divisible and is sent to zero by $H_{0}(-)$.

To summarize:

\begin{prop}\label{pr:HG}
Let $X$ be a spectrum and let $G$ be an abelian group.
A map $X \ra HG$ is phantom if and only if it is divisible
and is sent to zero by $H_{0}(-)$.    
The group of phantoms is isomorphic to $\PExt(H_{-1}X,G)$.   \qed
\end{prop}

\begin{proof}[Proof of Theorem~\ref{th:PExt}]
By the K\"unneth theorem, $H_{*}(V \Smash X)$ is a (non-canonical)
sum of tensor and torsion products of $H_{*}V$ and $H_{*}X$.
So the theorem follows from the following lemma.
\end{proof}

\begin{lemma}\label{le:PExt}
Let $A$, $B$ and $C$ be abelian groups with $C$ finitely generated.
If $\PExt(A,B) = 0$, then $\PExt(A \tensor C, B) = 0$ and 
$\PExt(A * C, B) = 0$, where $*$ denotes the torsion product.
\end{lemma}

\begin{proof}
We can assume without loss of generality that $C$ is cyclic.  
If it is infinite cyclic, then $A \tensor C = A$ and $A * C = 0$, 
so we are done.  
If it is $\Zn$, then both $A \tensor C$ and $A * C$ are killed by $n$, 
and so the $\Ext$ groups have no divisible elements.
\end{proof}

We now consider the finite type case.  

\begin{thm}\label{th:no}
Let $X$ and $Y$ be spectra such that there are
no phantom maps of any degree from $X$ to $Y$.  Let $V$ be a finite
spectrum.  Then $\Ph (V\Smash X,Y)_{*}$ is a bounded torsion graded
group.  In particular, if $Y$ has finite type, then there are no 
phantom maps from $V\Smash X$ to $Y$.
\end{thm}

\begin{proof}
For any finite spectrum $V$, there is a cofiber sequence 
\[
\bigvee _{i} S^{\alpha _{i}} \ra V \ra T ,
\]
where the wedge of spheres is finite and $T$ is a finite torsion spectrum.
In particular, $nT=0$ for some $n$.  Smashing this cofiber sequence with
$X$, we get a cofiber sequence
\[
\bigvee _{i} \Sigma ^{\alpha _{i}} X  \ra X \Smash V \ra X \Smash T .
\]
Our hypothesis implies that
$\Ph(\Wedge \Sigma^{\alpha_{i}} X,Y) = 0$.
Hence any phantom map $f\mathcolon X \Smash V \ra Y$ factors through a
map $g \mathcolon X \Smash T \ra Y$.  Since $ng=0$, it follows that
$nf=0$, so $\Ph (X \Smash V,Y)$ is a bounded torsion group, as
required. 

The last sentence of the theorem then follows immediately from
Theorem~\ref{th:div}. 
\end{proof}

Note that the second half of the theorem only relies on the divisibility
of $\Ph (V \Smash X,Y)$.   

\section{Chromatic phantom maps}\label{se:n-ph}

In this section we consider a generalization of phantom maps.  It is
much simpler to work $p$-locally for this generalization, where $p$ is
some fixed prime, so we will do so.  The symbol $\hat{X}$ will therefore
refer to $L_{M(p)}X$, the $p$-adic completion of $X$, rather than the
pro-finite completion of $X$.  Similarly, if $A$ is an abelian group,
$\hat{A}$ will refer to its $p$-completion.

Recall that a finite spectrum $V$ is said to be of \mdfn{type $n$} if
$K(n-1)_{*}V=0$ but $K(n)_{*}(V)$ is nonzero.  Let $\cat{C}_{n}$ denote
the thick subcategory of spectra of type at least $n$.  Then the
$\cat{C}_{n}$ form the complete list of nonzero thick subcategories of
finite spectra~\cite{hosm:nsht2}.

\begin{defn}
A map $X \ra Y$ is said to be \mdfn{$n$-phantom} if for every finite
spectrum $V\in \cat{C}_{n}$ and every map $V \ra X$ the composite $V \ra
X \ra Y$ is zero.  We write $\nPh (X,Y)$ for the group of $n$-phantom
maps from $X$ to $Y$.  
\end{defn}

Note that a $0$-phantom map is the same thing as a phantom map.  In
general, every $n$-phantom map is also an $n+1$-phantom map.  The
$n$-phantom maps form a two-sided ideal of morphisms, so that if $f$ is
$n$-phantom, then $gf$ and $fh$ are also $n$-phantom whenever we can
form the compositions.  Unlike ordinary phantom maps, this ideal is not
nilpotent when $n\geq 1$.  Indeed, if there are no maps at all from a
type $n$ spectrum to $X$, then the identity map of $X$ is $n$-phantom.
For example, the identity map of $K(n-1)$ is $n$-phantom but not
$n-1$-phantom.  This also shows that nonzero homotopy classes can
be $n$-phantom if $n>0$, though there are no nonzero $n$-phantoms whose
domain is a type $n$ finite spectrum, by definition.

Just as with phantom maps, there is a universal $n$-phantom map out of
$X$.  To see this, let $\cat{C}_{n}'$ denote a countable skeleton of
$\cat{C}_{n}$.  Let $P$ denote the coproduct over all maps
$V\xrightarrow{}X$ of $V$, where $V$ runs through spectra in
$\cat{C}_{n}'$.  Let $f\mathcolon X \ra Y$ denote the cofiber of the
obvious map 
$P \ra X$.  Then $f$ is obviously $n$-phantom, and if $g \mathcolon X
\ra Z$ is 
another $n$-phantom map, then $g$ factors (non-uniquely) through $f$.  In
particular, there are no $n$-phantom maps out of $X$ if and only if $X$
is a retract of a wedge of finite spectra of type at least $n$.  
For example, there is a nonzero $n$-phantom map out of any
finite spectrum of type $n-1$ or less.

There is also a universal $n$-phantom into $Y$, which will be discussed
in the next part of this section.

\subsection{Brown-Comenetz duality and universal $n$-phantom maps}

To characterize $n$-phantom maps, we must recall the finite localization
functor $L_{n-1}^{f}$ on the stable homotopy category.  
See~\cite{mi:fl} and~\cite[Section~3.3]{hopast:ash} for details.  
Recall that $L_{n-1}^{f}$ is a smashing localization whose acyclics
consist of the localizing subcategory generated by $\cat{C}_{n}$
and whose local objects are those spectra $Y$ such that $[V,Y] = 0$
for every $V$ in $\cat{C}_{n}$.
The fiber $C_{n-1}^{f}X$ of the map $X \ra L_{n-1}^{f}X$ is the closest 
one can get to $X$ using only finite spectra of type at least $n$.  
More precisely, let
$\cat{C}_{n}'$ denote a countable skeleton of $\cat{C}_{n}$, and let
$\Lambda _{n}X$ denote the category whose objects are maps $F \ra X$
where $F\in \cat{C}_{n}'$, and whose morphisms are commutative
triangles.  Then $\Lambda _{n}X$ is a filtered category, and the minimal
weak colimit of the obvious functor from $\Lambda _{n}X$ to finite
spectra is $C_{n-1}^{f}X$~\cite[Remark 2.3.18]{hopast:ash}.  
As a special case, $L_{-1}^{f} X = 0$ and $C_{-1}^{f} X = X$.

Since there are no maps from a finite $V$ in $\cat{C}_{n}$ to an
$L_{n-1}^{f}$-local spectrum $Y$, the identity map of such a $Y$ is
$n$-phantom.  More generally, we have the following proposition.  

\begin{prop}\label{pr:char-n}
A map $g \mathcolon X \ra Y$ is $n$-phantom if and only if $C_{n-1}^{f}g
\mathcolon 
C_{n-1}^{f}X \ra C_{n-1}^{f}Y$ is phantom.  
\end{prop}

\begin{proof}
Suppose first that $C_{n-1}^{f}g$ is phantom.  Suppose $V\in
\cat{C}_{n}$, and $h \mathcolon V \ra X$ is a map.  The composite $V \ra X \ra
L_{n-1}^{f}X$ is trivial, since $V$ is $L_{n-1}^{f}$-acyclic.  Hence $h$
factors through a map $V \ra C_{n-1}^{f} X$.  Since $C_{n-1}^{f} g$ is
phantom, the composite $V \ra C_{n-1}^{f} X \ra C_{n-1}^{f} Y$ is
trivial.  Hence the composite $V \ra X \ra Y$ is trivial, so $g$ is
$n$-phantom.  

Conversely, suppose $g$ is $n$-phantom, $V$ is finite, and $V \ra
C_{n-1}^{f}X$ is a map.  Since $C_{n-1}^{f}X $ is the minimal weak
colimit of the finite spectra in $\cat{C}_{n}$ mapping to $X$, the map
$V \ra C_{n-1}^{f}X$ must factor through a map $Z \ra C_{n-1}^{f}X$
where $Z$ is in $\cat{C}_{n}$.  Thus, we may as well assume that $V\in
\cat{C}_{n}$.  Then, since $g$ is $n$-phantom, the composite $V \ra
C_{n-1}^{f}X \ra C_{n-1}^{f} Y$ factors through a map $V \ra \Sigma^{-1}
L_{n-1}^{f}Y $.  But $V$ is $L_{n-1}^{f}$-acyclic, so this means the
composite $V \ra C_{n-1}^{f}X \ra C_{n-1}^{f} Y$ is zero.
\end{proof}

\begin{cor}\label{co:smash}
If $g$ is an $n$-phantom map, so is $g \Smash Z$ for all spectra $Z$.  
\end{cor}

\begin{proof}
Since $L_{n-1}^{f}$ is a smashing localization, we have
$C_{n-1}^{f}(g\Smash Z)=C_{n-1}^{f}(g)\Smash Z$.  The result now follows
from the fact that ordinary phantoms form an ideal under the smash
product, which is easily proved using Proposition~\ref{pr:char}.  
\end{proof}

The fact that ordinary phantom maps compose to 
$0$~\cite[Corollary~4.7]{chst:pmht} then yields the following corollary.  

\begin{cor}
Suppose $f \mathcolon X \ra Y $ and $g \mathcolon Y \ra Z$ are
$n$-phantom maps.  Then the 
composite $gf$ factors through a map $L_{n-1}^{f}X \ra Z$.  
\end{cor}

There is another localization functor which is relevant to $n$-phantom maps.
Let $F(n)$ denote a finite spectrum of type $n$.  
We can consider the localization functor $L_{F(n)}$ with respect 
to $F(n)$ and the associated acyclization functor $C_{F(n)}$,
both of which depend only on $n$.
These are orthogonal to $L_{n-1}^{f}$ and $C_{n-1}^{f}$ in the sense that 
the $F(n)$-acyclic spectra are the same as the $L_{n-1}^{f}$-local spectra.
(Such orthogonality is analyzed in~\cite[Theorem~3.3.5]{hopast:ash}.)
In particular, $L_{n-1}^{f}X$ is $F(n)$-acyclic for all $X$, 
and $C_{F(n)}X$ is $L_{n-1}^{f}$-local for all $X$.  
The spectrum $L_{F(n)}X$ should be thought of as the 
$(p,v_1,\dots,v_{n-1})$-completion of $X$.  
If $X$ is connective, then $L_{F(n)}X=L_{M(p)}X$~\cite[Theorem~3.1]{bo:lsh}.

In particular, the preceding corollary implies that if $f \mathcolon X
\ra Y$ and $g \mathcolon Y \ra Z$ are $n$-phantom maps, and either $Z$
is $F(n)$-local or $X$ is $L_{n-1}^{f}$-acyclic, then the composite $gf$
is null.  

\begin{prop}\label{pr:no-n}
There are no nonzero $n$-phantom maps to $Y$ if and only if there are no
nonzero phantom maps to $Y$ and $Y$ is $F(n)$-local.  
\end{prop}

\begin{proof}
Suppose first that there are no nonzero $n$-phantom maps to $Y$.  Then
obviously there are no nonzero phantom maps to $Y$.  Also, the identity
map of $C_{F(n)}Y$ is $n$-phantom, since $C_{F(n)}Y$ is
$L_{n-1}^{f}$-local.  Hence the map $C_{F(n)}Y \ra Y$ must be $0$, so
$Y$ is a retract of $L_{F(n)}Y$.  This implies that $Y$ is
$F(n)$-local.  

Conversely, suppose there are no nonzero phantom maps to $Y$ and $Y$ is
$F(n)$-local.  Suppose $X \ra Y$ is an $n$-phantom map.  Then the
composite $C_{n-1}^{f}X \ra X \ra Y$ is phantom by
Proposition~\ref{pr:char-n}.  Hence it is $0$, so our original map
factors through a map $L_{n-1}^{f}X \ra Y$.  But $L_{n-1}^{f}X$ is
$F(n)$-acyclic, so this map must be $0$.  
\end{proof}

Proposition~\ref{pr:no-n} suggests the following definition, which is
due to Mahowald and Rezk~\cite{mare}.  

\begin{defn}
Given $n\geq 0$, define a functor $W_{n-1}$ on \ulp $p$-local\urp\ spectra by
$W_{n-1}Y=IC_{n-1}^{f}Y = \FF(Y, W_{n-1}S^{0})$.  
\end{defn}

Note that $W_{-1}Y=IY$, and that $W_{n-1}Y=0$ if and only if $Y$ is
$F(n)$-acyclic, or equivalently, $L_{n-1}^{f}$-local.  Note as well that
the natural map $C_{n-1}^{f}Y \ra Y$ gives rise to a natural map $IY \ra
W_{n-1}Y$.

\begin{lemma}\label{le:W-F(n)}
The natural map $IY \ra W_{n-1}Y$ is the $F(n)$-localization of $IY$.  
\end{lemma}

\begin{proof}
The fiber of the natural map $IY \ra W_{n-1}Y$ is $IL_{n-1}^{f}Y$.  This
is an $L_{n-1}^{f}S^{0}$-module spectrum, so is $F(n)$-acyclic.  
It therefore suffices to show that $W_{n-1}Y$ is $F(n)$-local.  Suppose $Z$
is $F(n)$-acyclic.  Then $[Z,W_{n-1}Y]=[Z\Smash C_{n-1}^{f}Y,I]=0$.
Indeed, $C_{n-1}^{f}Y$ is built out of finite spectra of type at least
$n$, so since $Z \Smash F(n)=0$, we must have $Z\Smash C_{n-1}^{f}Y=0$.
\end{proof}

\begin{cor}\label{co:BC-n}
There are no nonzero $n$-phantom maps to $W_{n-1}Y$.  
\end{cor}

\begin{proof}
We have $W_{n-1}Y=IC_{n-1}^{f}Y=L_{F(n)}IY$, so there are no nonzero
phantom maps to $W_{n-1}Y$ and $W_{n-1}Y$ is $F(n)$-local.
Proposition~\ref{pr:no-n} completes the proof.  
\end{proof}

It follows that there are no nonzero $n$-phantom maps from
$X$ to $IY$ if $L_{n-1}^{f}(X \Smash Y)=0$.
For if $X \ra IY$ is an $n$-phantom map, then so is its
adjoint $X \Smash Y \ra I = IS^{0}$, 
as it factors through $X \Smash Y \ra IY \Smash Y$.
By Corollary~\ref{co:BC-n}, the map $X \Smash Y \ra I$ factors 
through $C_{F(n)}I$, which is $L_{n-1}^{f}$-local.  
But $X \Smash Y$ is $L_{n-1}^{f}$-acyclic, so the map is zero.

As a special case of the Corollary, 
there are no nonzero $n$-phantom maps to $IV=W_{n-1}V$
when $V$ is a finite spectrum of type at least $n$, since
$C_{n-1}^{f}V=V$.

\begin{prop}\label{pr:universal}
There is a universal $n$-phantom map to $Y$.  
\end{prop}

\begin{proof}
Define $Q$ to be the product over all maps $Y \ra IV$ of $IV$,
where $V$ runs through spectra in $\cat{C}_{n}'$.  Let $g \mathcolon F
\ra Y$ 
denote the fiber of the obvious map $Y \ra Q$.  Then any $n$-phantom map
to $Y$ will have to factor through $g$, since there are no $n$-phantom
maps to $Q$.  We claim that $g$ is itself $n$-phantom.  Indeed, suppose
$V \ra Y$ is a nonzero map, where $V$ is a finite spectrum of type at
least $n$.  This map corresponds to a homotopy class in $Y\Smash DV$,
and so there is some map $Y\Smash DV \ra I$ which does not send this
homotopy class to $0$.  The adjoint of this map is a map $Y \ra IDV$,
and the composite $V \ra Y \ra IDV$ is nonzero, since its adjoint $V
\Smash DV \ra I$ sends the unit to a nontrivial class in $\Q /\Z $.
Thus, if $V \ra Y$ is nonzero, so is the composite $V \ra Y \ra Q$, and
so $g$ is $n$-phantom.
\end{proof}

This construction shows that there are no $n$-phantom maps to $Y$ if and
only if $Y$ is a retract of a product of Brown-Comenetz duals of finite
spectra in $\cat{C}_{n}$.  

Just as in the usual case, there is a natural map $Y \ra W_{n-1}^{2}Y$
adjoint to the evaluation map $Y \Smash W_{n-1}Y \cong Y \Smash
\FF(Y,W_{n-1}S^{0}) \ra W_{n-1}S^{0}$.  Note that this map is the
composite $Y \ra I^{2}Y \ra L_{F(n)}I^{2}Y$.  Indeed, we have
$W_{n-1}^{2}Y = L_{F(n)}I^{2}C_{n-1}^{f}Y$.  But $I^{2}L_{n-1}^{f}Y$ is
$F(n)$-acyclic, so $L_{F(n)}I^{2}C_{n-1}^{f}Y \cong L_{F(n)}I^{2}Y$.   

\begin{prop}\label{pr:double-n}
The fiber of the natural map $Y \ra W_{n-1}^{2}Y = L_{F(n)}I^{2}Y$ is a
universal $n$-phantom map into $Y$.
\end{prop}

\begin{proof}
Corollary~\ref{co:BC-n} implies that there are no nonzero $n$-phantom
maps into $W_{n-1}^{2}Y$.  To complete the proof, it therefore
suffices to show that the map $[V,Y] \ra [V, W_{n-1}^{2}Y]$ is
injective for all finite spectra $V$ of type at least $n$.  It is easy
to check that $[V,Y] \ra [V,I^{2}Y]$ is injective for all finite
spectra.  But $[V,I^{2}Y] \ra [V,L_{F(n)}I^{2}Y]$ is an isomorphism for
all finite spectra of type at least $n$.  
\end{proof}

An interesting consequence of Proposition~\ref{pr:double-n} is that, if
there are no nonzero phantom maps to $Y$, then there are no nonzero
$n$-phantom maps to $L_{F(n)}Y$.  Indeed, if there are no nonzero
phantom maps to $Y$, then $Y \ra I^{2}Y$ is a split monomorphism.  It
follows that $L_{F(n)}Y \ra L_{F(n)}I^{2}Y$ is a split monomorphism, and
so, by Proposition~\ref{pr:double-n}, there are no nonzero $n$-phantom
maps to $L_{F(n)}Y$.  

These results suggest that $W_{n-1}$ is the right notion of
Brown-Comenetz duality in the $F(n)$-local stable homotopy category.
Recall from~\cite[Appendix~B]{host:mktl} that this category is an
algebraic stable homotopy category in its own right, though the sphere
$L_{F(n)}S^{0}=\hat{S^{0}}$ is not small.  Any finite spectrum of type
$n$ is a small weak generator.  The natural notion of phantom map in
this setting is therefore the notion of an $n$-phantom map.  The
comments preceding Proposition~\ref{pr:no-n} imply that the composition
of two $n$-phantoms in this local category is null, which also follows
from~\cite[Theorem~4.2.5]{hopast:ash}.  Ordinary Brown-Comenetz duality
will not stay in this local subcategory, but $W_{n-1}$ obviously does do
so.

Note that the Brown-Comenetz duality $\hat{I}$ in the $K(n)$-local
category studied in~\cite[Section~10]{host:mktl} is just the restriction
of $W_{n-1}$.  Indeed, by definition, $\hat{I}Y=IM_{n}Y$ for a
$K(n)$-local spectrum $Y$.  But, for an $L_{n}$-local spectrum $Y$,
$M_{n}Y=C_{n-1}^{f}Y$, so $\hat{I}Y=W_{n-1}Y$.    

We now introduce one possible analog of the finite type condition,
appropriate for chromatic phantom maps.

\begin{defn}
Define a spectrum $Y$ to be of \mdfn{$n$-finite type} if $Y\Smash V$ 
has finite homotopy groups for all $V$ of type at least $n$.
\end{defn}

Note that the collection of all $V$ such that $Y\Smash V$ has finite
homotopy groups is a thick subcategory, so that $Y$ is of $n$-finite type if 
and only if there is a type $n$ finite spectrum $V$ such that $Y\Smash V$ 
has finite homotopy groups.

If $Y$ is of $n$-finite type, then $Y$ is of $(n+1)$-finite
type.  The Johnson-Wilson spectrum $E(n)$, whose homotopy groups are
given by 
\[
E(n)_{*}\cong \Z _{(p)}[v_{1},\dots ,v_{n-1}][v_{n},v_{n}^{-1}] ,
\]
is of $n$-finite type but not of $(n-1)$-finite type.  

We also have the following lemma, whose validity is the main reason for
introducing spectra of $n$-finite type.  

\begin{lemma}\label{le:lns0}
The spectrum $L_{n}S^{0}$ is of $n$-finite type for $n\geq 1$.  
\end{lemma}

\begin{proof}
Suppose $V$ is a type $n$ finite spectrum.  Then $L_{n}S^{0}\Smash
V=L_{n}V=L_{K(n)}V$.  Corollary~8.12 of~\cite{host:mktl} shows that
$\pi _{k}L_{K(n)}V$ is finite for all $k$.  
\end{proof}

It follows that $L_{n}Y$ is of $n$-finite type for any finite spectrum
$Y$ and any $n\geq 1$.  Any dualizable spectrum in the $K(n)$-local
category is of $n$-finite type by~\cite[Theorem~8.9]{host:mktl}.  
A connective $p$-complete spectrum is
of $n$-finite type if and only if it is an fp-spectrum of fp-type at
most $n-1$, in the sense of~\cite{mare}.

In contrast to Lemma~\ref{le:lns0}, 
the calculations of~\cite{marash:vphctc} suggest that
$L_{2}^{f}S^{0}$ is not of $2$-finite type.

It is widely believed that in fact $L_{n}S^{0}$ is of $1$-finite type
for all $n$.  This is true for $L_{2}S^{0}$ if $p\geq 5$, by the results
of Shimomoura~\cite{sh:ans} as explained in~\cite[Theorem~15.1]{host:mktl}.  
This conjecture is implied by~\cite[Conjecture~3.8]{mare}, and by the
chromatic splitting conjecture of~\cite{ho:blfhcsc}.  

For us, the main advantage of spectra of $n$-finite type is the
following theorem.  

\begin{thm}\label{th:n-ph}
Suppose $Y$ is of $n$-finite type.  Then the map $Y \ra
W_{n-1}^{2}Y$ is $F(n)$-localization.  In particular, there are no
$n$-phantom maps to $L_{F(n)}Y$, and the inclusion $C_{F(n)}Y \ra Y$ is
a universal $n$-phantom map to $Y$.
\end{thm}

\begin{proof}
Consider the map $Y \ra W_{n-1}^{2}Y=L_{F(n)}I^{2}Y$.
The target of this map is obviously $F(n)$-local.  On the other hand, if
we smash with $F(n)$, we have $L_{F(n)}I^{2}Y \Smash
F(n)=L_{F(n)}I^{2}(Y\Smash F(n))=Y\Smash F(n)$, since $Y$ is of
$n$-finite type.
Thus $Y \ra W_{n-1}^{2}Y$ is $F(n)$-localization.

The remaining statements follow from Proposition~\ref{pr:double-n}.
\end{proof}

A corollary of this theorem is that $W_{n-1}$ defines a contravariant
self-equivalence of the category of $F(n)$-local spectra of $n$-finite
type.  
This generalizes~\cite[Corollary~8.3]{mare}, since for $n \geq 1$
a connective spectrum is $F(n)$-local if and only if it is $p$-complete.

If we apply Theorem~\ref{th:n-ph} to $L_{n}S^{0}$ we get the following
corollary.  

\begin{cor}\label{co:Ln}
The kernel of the map $\pi _{*}L_{n}S^{0}\ra \pi _{*}L_{K(n)}S^{0}$ is
the subgroup of $n$-phantom homotopy classes.  Furthermore, there are no
$n$-phantom maps to $L_{K(n)}S^{0}$.  
\end{cor}

Similarly, there are no $n$-phantom
maps to $L_{K(n)}E(n) = L_{F(n)}E(n)$, and the kernel of the map 
$E(n)^{*}X \ra (L_{K(n)}E(n))^{*}(X)$ is the subgroup of $n$-phantom 
cohomology classes.  

If $L_{n}S^{0}$ has $1$-finite type, then the $k$-phantom homotopy
classes in $\pi _{*}L_{n}S^{0}$ form the kernel of the map 
\[
\pi _{*}L_{n}S^{0} \ra \pi _{*}L_{F(k)}L_{n}S^{0}=\pi_* L_{K(k)\vee \dots \vee
K(n)}S^{0} .
\]
The chromatic splitting conjecture~\cite{ho:blfhcsc}, if true, implies
that this kernel is built from $\pi _{*}L_{j}S^{0}$ for $j<k$ in a
predictable way.  In particular, the $1$-phantom subgroup of 
$\pi_{*}L_{n}S^{0}$ should be the $\Q /\Z _{(p)}$ summands ($\Q$ for $n = 0$),
and the chromatic splitting conjecture would tell us that there are $2^{n}-1$ 
of these located in specific dimensions.

\subsection{Divisibility}

We now describe the notion of divisibility which is relevant to $n$-phantom
maps.

\begin{defn}\label{de:vn-div}
A map $X \ra Y$ is \mdfn{$v_{n-1}$-divisible} if for every finite 
spectrum $V$ of type at least $n-1$, every $v_{n-1}$-self map 
$v \mc V \ra \Sigma^{-d} V$, and every map $V \ra X$, the composite 
$V \ra Y$ is $v$-divisible.
By \mdfn{$v$-divisible} we mean that the map factors through
$v^{k} \mc V \ra \Sigma^{-dk} V$ for each $k$.
\end{defn}

The centrality of $v_{n-1}$-self maps~\cite{hosm:nsht2} implies that
testing divisibility against any single $v_{n-1}$-self map of $V$
ensures divisibility with respect to any other $v_{n-1}$-self map of $V$.

\begin{prop}\label{pr:vn-div}
A map $X \ra Y$ is $v_{n-1}$-divisible if and only if 
it is $n$-phantom.
\end{prop}

\begin{proof}
If $X \ra Y$ is $v_{n-1}$-divisible and $V \ra X$ is a map from a
finite spectrum of type at least $n$, then the zero map is a 
$v_{n-1}$-self map of $V$.  Therefore the composite $V \ra Y$ is zero.

Conversely, suppose that $X \ra Y$ is $n$-phantom and let $V \ra X$
be a map from a finite spectrum of type at least $n-1$. 
Let $v \mc V \ra \Sigma^{-d} V$ be a $v_{n-1}$-self map of $V$
and write $V/v^{k}$ for the cofibre of $v^{k}$.
This cofibre has type at least $n$, and so the composite
$\Sigma^{-1} V/v^{k} \ra V \ra X \ra Y$ is trivial.
Thus $V \ra X \ra Y$ is divisible by $v^{k}$.
\end{proof}

A divisible map is clearly $v_{0}$-divisible, so we obtain the following
corollary.

\begin{cor}
A divisible map is $1$-phantom.
\end{cor}

Since there are non-divisible $1$-phantoms, there are non-divisible maps
which are $v_{0}$-divisible.

Sometimes, the converse of the Corollary holds.

\begin{prop}\label{pr:one}
If $X \Smash M(p)$ is finite or $Y$ has $1$-finite type, 
then a map $X \ra Y$ is divisible if and only if it is $v_{0}$-divisible
if and only if it is phantom.
In case $Y$ has $1$-finite type, $\oPh (X,Y)$ is
the divisible subgroup of $[X,Y]$.  
\end{prop}

\begin{proof}
We have already observed that divisible implies $v_{0}$-divisible
and that $v_{0}$-divisible is equivalent to being $1$-phantom.

Suppose that  $X \ra Y$ is a $1$-phantom map.
If $X \Smash M(p)$ is finite, then so is $X\Smash M(p^{n})$ for all $n$.
Hence the composite $X \Smash \Sigma^{-1}M(p^{n}) \ra X \ra Y$ must be
$0$ for all $n$, so the map $X \ra Y$ must be divisible by $p^{n}$ for
all $n$.  On the other hand, if $Y$ has $1$-finite type, then there are
no $1$-phantom maps to $L_{M(p)}Y$.  Hence any $1$-phantom map factors
through the rational spectrum $C_{M(p)}Y$, so is divisible as a
$1$-phantom map.  
\end{proof}

For example, take $Y=L_{1}S^{0}=L_{1}^{f}S^{0}$~\cite[Theorem~10.12]{ra:lrc}.  
Then all the homotopy groups of $L_{1}S^{0}$ are finitely generated,
except $\pi_{-2}$, which has a $\Q/\Z_{(p)}$ 
summand~\cite[Theorem~8.10]{ra:lrc}.  
This summand is the $1$-phantom subgroup.  

Here is a variant of Proposition~\ref{pr:vn-div}.

\begin{prop}
Let $V$ be a finite spectrum of type at least $n-1$, where $n \geq 1$,
and let $v \mathcolon V \ra \Sigma^{-d} V$ be a $v_{n-1}$-self map.  
Then a map $f \mathcolon V \ra Y$ is $n$-phantom 
if and only if it is $v$-divisible.  Furthermore, if
$Y$ has $n$-finite type, then the group of $n$-phantom maps from $V$ to
$Y$ is itself $v$-divisible.
\end{prop}

\begin{proof}
By the centrality of $v_{n-1}$-self maps~\cite{hosm:nsht2},
$f$ is $v$-divisible if and only if it is $v_{n-1}$-divisible,
and so the first part follows from Proposition~\ref{pr:vn-div}.

Now suppose that $Y$ has $n$-finite type and $f \mathcolon V \ra Y$ is
$n$-phantom.  A map to $Y$ is $n$-phantom if and only if it factors
through $C_{F(n)}Y$, which is $L_{n-1}^{f}$-local.  Hence $f$ factors
through $L_{n-1}^{f}V=v^{-1}V$.
Since $v$ is a self-equivalence of
$v^{-1}V$, $f$ is $v$-divisible as a map into $v^{-1}V$, and hence as a
map into $C_{F(n)}V$, and hence as an $n$-phantom map.
\end{proof}

Now we introduce another finite type condition.

\begin{defn}\label{de:vn-ft}
A spectrum $Y$ has \mdfn{$v_{n}$-finite type} if for each finite
spectrum $V$ of type $n$ and each $v_{n}$-self map $v \mc V \ra \Sigma^{-d} V$,
the graded group $[V,Y]_{*}$ contains no nonzero $v$-divisible elements.
\end{defn}

Testing this condition against any single $v_{n}$-self map of $V$
ensures that it holds with respect to all other $v_{n}$-self maps of $V$.

A spectrum of finite type 
is always of $v_{0}$-finite type.  However, the converse is false,
as one can see by taking a large product of suitable finite type
spectra.  In the same way, one sees that $v_{n}$-finite type does
not imply $n$-finite type.
For $n \geq 1$ it is also the case that $n$-finite type does not
imply $v_{n}$-finite type.
Indeed, $K(n)$ has $n$-finite type but not $v_{n}$-finite type.
Similar examples show that when $m \neq n$, $v_{m}$-finite type
and $v_{n}$-finite type are unrelated.

It is clear that products and coproducts of spectra of 
$v_{n}$-finite type also have $v_{n}$-finite type.
We will see below that the sphere has $v_{n}$-finite type for each $n$.
Since there exist spectra which do not have $v_{n}$-finite type, it
follows that for each $n$ the collection of spectra of $v_{n}$-finite
type is not thick.
For an explicit example which shows that the collection of $v_{0}$-finite
type spectra is not thick, chose a projective resolution of $\Q$ and
look at the corresponding cofibre sequence of \EM spectra.

Our reason for introducing the $v_{n}$-finite type condition is the
next result.

\begin{prop}\label{pr:vn-ft}
The following conditions on a spectrum $Y$ are equivalent:
\begin{roenumerate}
\item $Y$ has $v_{n}$-finite type.
\item For every $X$ and every map $f \mc X \ra Y$, 
      $f$ is $n$-phantom if and only if it is $n+1$-phantom.
\end{roenumerate}
\end{prop}

\begin{proof}
(i) $\implies$ (ii):
First recall that every $n$-phantom is $n+1$-phantom.  
So assume that $X \ra Y$ is $n+1$-phantom and let $V \ra X$ be a map 
from a type $n$ finite spectrum with a $v_{n}$-self map $v$.
By Proposition~\ref{pr:vn-div}, the composite $V \ra Y$ is
$v$-divisible.
Since $Y$ has $v_{n}$-finite type, this composite is trivial.
Thus the original map is $n$-phantom.

(ii) $\implies$ (i):  Assume that a map $X \ra Y$ is $n$-phantom
if and only if it is $n+1$-phantom.  We will prove that $Y$ has
$v_{n}$-finite type.  Let $V$ be a type $n$ finite spectrum with
a $v_{n}$-self map $v$, and let $f \mc V \ra Y$ be a $v$-divisible
map.  We must show that $f$ is trivial.
By the centrality of $v_{n}$-self maps, $f$ is $v_{n}$-divisible.
By Proposition~\ref{pr:vn-div}, this implies that $f$ is $n+1$-phantom.
By assumption, this implies that $f$ is $n$-phantom.
But $V$ is type $n$, so $f$ must be trivial.
\end{proof}

There are examples of spectra of $v_{n}$-finite type.
In order to describe them, we let $\cD$ denote the class of abelian
groups which are the sum of a finitely generated free group and a
bounded torsion group.
One can show that $\cD$ is closed under kernels, cokernels
and extensions.

\begin{prop}\label{pr:n=1=0}
\begin{roenumerate}
\item If $Y$ is bounded below, then $Y$ has $v_{n}$-finite type for
each $n \geq 1$ and a map $X \ra Y$ is $n$-phantom if and only if it
is $1$-phantom.
\item If $Y$ is such that, for each $k$, $\pi_{k} Y$ is in $\cD$,
then $Y$ has $v_{0}$-finite type and a map $X \ra Y$ is $1$-phantom if
and only if it is phantom.
\end{roenumerate}
Combining the above, we see that if $Y$ is a bounded below spectrum
such that each $\pi_{k} Y$ is in $\cD$, then $Y$ has $v_{n}$-finite
type for each $n$ and a map $X \ra Y$ is $n$-phantom if and only if
it is phantom.
\end{prop}

\begin{proof}
(i) If $Y$ is bounded below and $V$ is finite, then for large enough
$m$ there are no nonzero maps from $\Sigma^{-m} V$ to $Y$.
Thus $Y$ has $v_{n}$-finite type for $n \geq 1$.
That $n$-phantoms are $1$-phantoms follows from Proposition~\ref{pr:vn-ft}.

(ii) If each $\pi_{k} Y$ is in $\cD$, then for any finite $V$
the group $[V,Y]$ is in $\cD$ as well.  But the groups in $\cD$
have no divisible elements, so $Y$ has $v_{0}$-finite type.
That $1$-phantoms are phantoms follows from Proposition~\ref{pr:vn-ft}.
\end{proof}

Propositions~\ref{pr:one} and~\ref{pr:n=1=0} give another proof of
Theorem~\ref{th:div}.

Note that there exist spectra $Y$ of $v_{0}$-finite type for which
it is not the case that $\pi_{k} Y$ is in $\cD$.
For example, take $Y = \prod_{l} H\Zpl$.

Consider the case $Y=E(n)$.  Let $I$ denote the ideal $(p,v_{1},\dots
,v_{n-1})\in E(n)_{*}$.  Define a class in $E(n)^{*}(X)$ to be
\mdfn{$I$-divisible} if it is in $I^{\infty }E(n)^{*}(X)=\bigcap _{k}
I^{k}E(n)^{*}(X)$.  Then the $I$-divisible cohomology classes should be
the $n$-phantom cohomology classes.  We can not completely prove this,
but we do have the following proposition.  

\begin{prop}
If a map $X \ra E(n)$ is $I$-divisible, then it is $n$-phantom.
\end{prop}

\begin{proof}
Suppose that $f \mathcolon X \ra E(n)$ is $I$-divisible, and $g
\mathcolon V \ra X$ is a 
map where $V\in \cat{C}_{n}$. Then $gf$ is an $I$-divisible element of
$E(n)^{*}(V)$.  However, for a particular type $n$ spectrum of the form
$M=M(p^{i_{0}},v_{1}^{i_{1}}, \dots , v_{n-1}^{i_{n-1}})$, one can easily
see that $E(n)^{*}(M)$ is all killed by a power of $I$.  It follows by a
thick subcategory argument that $E(n)^{*}(V)$ is all killed by a power
of $I$ for any $V\in \cat{C}_{n}$.  Thus any $I$-divisible element in
$E(n)^{*}(V)$ is zero, and so $f$ is $n$-phantom.  
\end{proof}

We do not know if the converse to this proposition holds.  It is
certainly true that if $X \ra E(n)$ is $n$-phantom, then the composite
$X \ra E(n) \ra E(n) \Smash M(p^{i_0},v_1^{i_1},\dots,v_{n-1}^{i_{n-1}})$ 
is trivial, since the second map factors through $L_{F(n)}E(n)$
and there are no $n$-phantom maps to $L_{F(n)}E(n)$.  But it is not
clear that this implies that our map is $I$-divisible.  

Finally, if $Y$ is an \EM spectrum, we have the following
characterization of $n$-phantom maps, analogous to the results
of~\cite[Section~6]{chst:pmht}.

\begin{prop}\label{pr:EM}
Suppose $B$ is an Abelian group.  If $n\geq 1$, the $n$-phantom subgroup
of $[X,HB]$ is precisely the subgroup of divisible elements.
Furthermore, an $n$-phantom map $X \ra HB$ is phantom if and only if the
map $H_{0}X \ra B$ is zero.  
\end{prop}

\begin{proof}
Since $HB$ is bounded below, the $n$-phantom subgroup coincides with the
$1$-phantom subgroup, by Proposition~\ref{pr:n=1=0}.  Every divisible
map is $1$-phantom by Proposition~\ref{pr:one}.  So suppose $X \ra HB$
is $1$-phantom.  Then the composite $X \ra HB \ra HB \Smash M(p^{n})$ is
$1$-phantom, and hence phantom by Proposition~\ref{pr:n=1=0}.  But
$HB \Smash M(p^{n})=H(B/p^{n})\vee \Sigma HC$, where $C$ is the group
of $p^{n}$-torsion elements in $B$.  In particular, there are no phantom
maps to $HB \Smash M(p^{n})$.  Indeed, a phantom map from X to an \EM 
spectrum $HA$ is divisible, by Proposition~\ref{pr:HG}.
Since both $H(B/p^{n})$ and $HC$ are killed by $p^{n}$, 
they are not the target of a nonzero phantom map.
Hence the composite $X \ra HB \ra HB \Smash M(p^{n})$ is null
for all $n$, and so $X \ra HB$ is a divisible map.  

The last statement is just Proposition~\ref{pr:HG}.
\end{proof}

For example, if $X=HA$, then $\Ph(HA,HB)$ is $\PExt(A,B)$
concentrated in degree $1$, and $\oPh (HA,HB)$ is $\PExt(A,B)$ in degree
$1$ and the subgroup of divisible elements in $\Hom (A,B)$ in degree
$0$.

\renewcommand{\baselinestretch}{1}\normalsize
\bibliography{jourabbrev,christensen}
\bibliographystyle{hplain}

\end{document}